\documentclass[12pt]{article}
\usepackage[utf8]{inputenc}
\usepackage{latexsym}     
\usepackage{amsfonts}
\usepackage{amsthm}
\usepackage{amsmath}
\usepackage{amssymb}
\usepackage{graphicx} 
\newtheorem{theorem}{Theorem}[section]

\theoremstyle{remark}

\def\tm{{\pitchfork}}
\def\mt{{\emptyset}}

\def\Q{{\mathbb Q}}

\def\d{\partial}

\def\emp{{\emptyset}}

\def\Z{\mathbb{Z}}
\def\R{\mathbb{R}}
\def\N{\mathbb{N}}

\title{Infinite-order combinatorial Transverse Intersection Algebra TIA  via the probabilistic wiggling* model}
\author{Daniel An\footnote{SUNY Maritime College}, Ruth Lawrence\footnote{Einstein Institute of Mathematics, Hebrew University of Jerusalem},\  and Dennis~Sullivan\footnote{CUNY Graduate Center NY \& SUNY, Stony Brook NY}}
\date{12/03/25\footnote{This work was supported by the US-Israel Binational Science Foundation, Jerusalem}}

\begin{document}

\maketitle
\begin{abstract}
  \scriptsize  
  This paper constructs  a graded-commutative, associative, differential Transverse Intersection Algebra TIA {on the torus (in any dimension) with its cubical decomposition by using a probabilistic wiggling \footnote{Jiggling below refers to  translating back and forth in various directions, wiggling above and below  encompasses jiggling plus stretching and squeezing, these all weighted with  probability distributions.} interpretation.  This structure agrees with the combinatorial graded intersection algebra (graded by codimension) defined by transversality on pairs of `cuboidal chains'  which are in  general position. In order to define an intersection of cuboids which are not necessarily in general position, the boundaries of the cuboids are considered to be `wiggled' by a distance small compared with the lattice parameter, according to a suitable probability distribution and then almost always the wiggled cuboids will be in general position, producing a transverse intersection with new probability distributions on the bounding sides.  In order to make a closed theory, each geometric cuboid appears in an infinite number of forms with different probability distributions on the wiggled boundaries. The resulting structure is commutative, associative and satisfies the product rule with respect to the natural boundary operator deduced from the geometric boundary of the wiggled cuboids.

 This TIA can be viewed as a  combinatorial analogue of differential forms in which the continuity of space has been replaced by a lattice with corrections to infinite order. See  comparison to  Whitney forms  at the end of the paper and in \cite{S77}, \cite{Boss}, \cite{H}, \cite{A}.

 In order to obtain finite approximations it is possible at any given order to divide out by an ideal generated by elements of higher order than the given one and then Leibniz will only hold partially up to that order, due to the boundary operator not preserving the ideal.

  For application to fluid algebra we also consider the same construction starting with the $2h$ cubical complex instead of the $h$ cubical complex. The adjoined higher order elements will be identical to those required in the $h$ cubical complex.
  
  The $d$-dimensional theory is a tensor product of $d$ copies of the one-dimensional theory.
  }
\end{abstract}

\section{Aim and background}

The legendary paper of Ren\'e Thom \cite{Thom} lead to the phrase ``Thom transversality" as in Milnor's beautiful notes on that theory (\cite{Milnor}, also see Princeton Math Dept  notes of same). The main geometric idea of \cite{Thom} being Thom's continuously differentiable transversality  which was  $C^1$ generic, local and relative in contrast to the very attractive Sard theorem  which depended on higher smoothness needed for arbitrary dimensions. 
\

Geometric topology developed with these Thom 
 Transversality pictures behind every  geometric forehead in the  50's 60's 70's etc.
 ``Transversality"  was the key geometric word inside Thom's paper and is the motivating word for this paper. 

Remark: The above is true even though the competing word ``cobordism" created in  Thom's paper (by transversality) is even more legendary in algebraic and geometric topology, e.g. cobordism theory being 
 the first generalized homology theory and also the basis for the first proof of the Atiyah--Singer Index theorem (the index being an appropriate cobordism invariant of a geometrically defined  operator on a closed manifold). 

The full discrete analog of  such   geometric  and analytic  theories, like the transversality in this paper   and the exterior product of differential forms  with the analogy discussed in \cite{S77} is not yet forthcoming (but see \cite{Teleman}).

\

Let us begin the discretization of the dga of  differential forms or rather    
  this paper's dual (in the sense of Poincar\'e)  transversal intersection algebra  with differential satisfying the  Leibniz rule.

\

Consider a cubical lattice with lattice spacing $h$ in three dimensions. The usual associated chain complex has four non-trivial chain spaces, in dimensions 0,1,2 and 3 which have basis elements (see Figure 1) which are points, elemental edges of length $h$ parallel to one of the three axes, elemental plaquettes which are squares of edge length $h$ parallel to one of the three coordinate planes, and elemental cubes of edge length $h$, respectively. 
\[
\includegraphics[width=.6\textwidth]{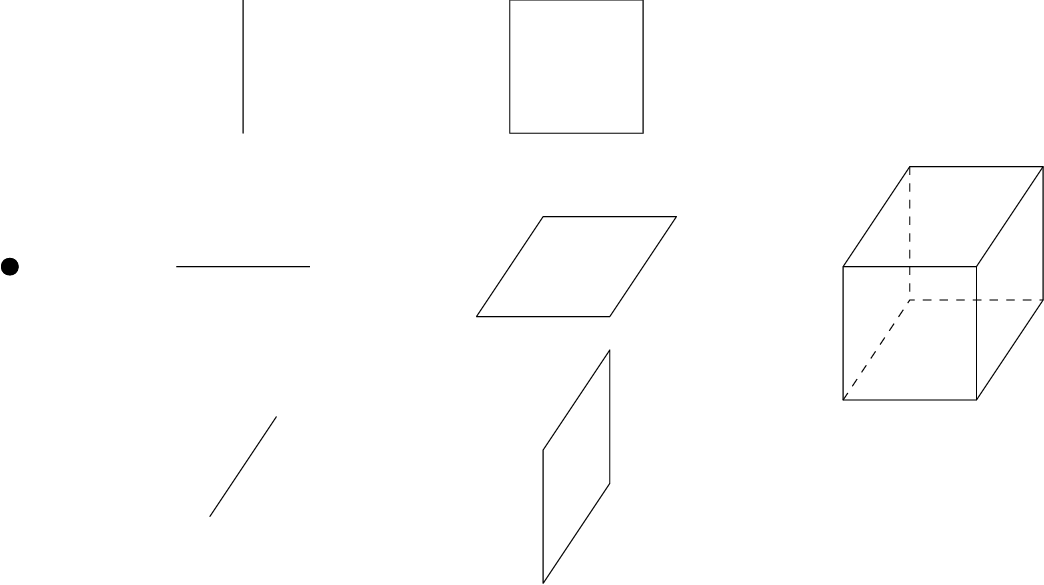}
\]\[\hbox{\small\sl Figure 1: The $h$ complex}
\]

\medskip
\noindent The boundary map is defined by the geometric boundary (with orientations).
 
 For application to fluid algebras we will also consider the $2h$ complex which has basis consisting of all possible similar cells but with edge length $2h$, as in Figure 2.
\[ 
\includegraphics[width=.5\textwidth]{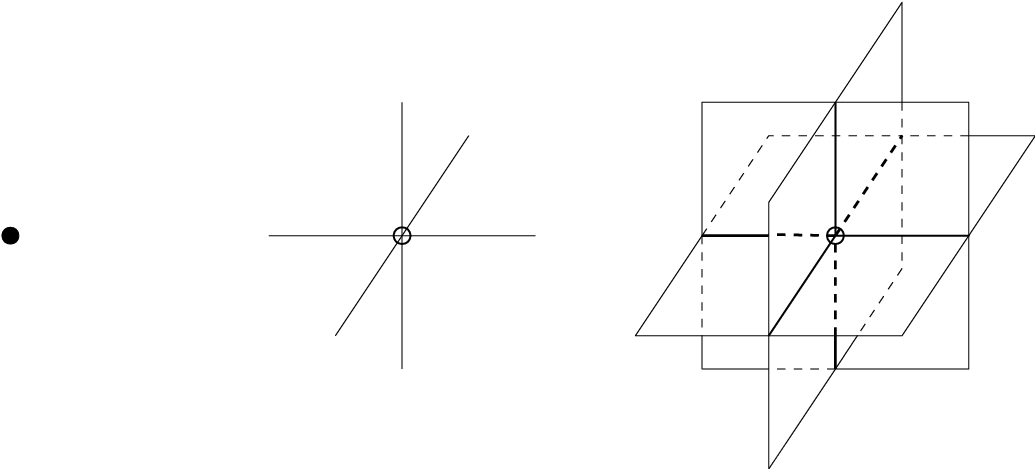}\hskip2em
\includegraphics[width=.22\textwidth]{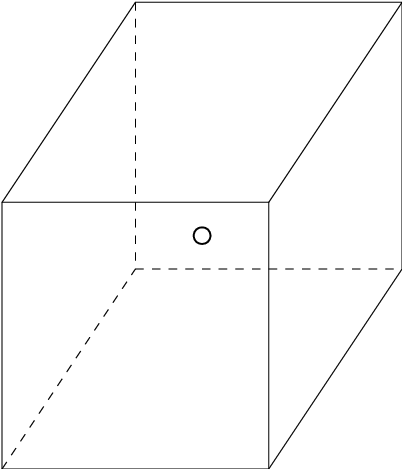}
\]
\[\hbox{\small\sl Figure 2: The $2h$ complex}\]

The $2h$ complex has nice properties:  each vertex is the barycenter of (1,3,3,1) $2h$ cells of dimensions 0,1,2,3 as in   Figure 2  and these form an  exterior algebra structure at each vertex. The entire $2h$  chain complex with boundary operator of degree $-1$  has a graded commutative associative (with the sign determined by the codimensions) intersection algebra structure  with   degrees $i,j$ going  to  degree $i+j -3$. Plus there is a natural star duality relating degrees one and two and  relating degrees zero and three. 

These two chain complexes, the $h$ cubical decomposition and the complex of overlapping $2h$ cells  give a first approximation to the  dual intersection  geometric picture of the  exterior algebra structure on  differential forms with the geometric  boundary operator of degree $-1$  being in  Hom duality with  the  picture of the exterior  derivative on forms \cite{S}. 

 The new caveat   is that the $\d$ operator of degree $-1$   does not satisfy  the product rule for this first approximation of the geometric product.  This  discrepancy   has been  treated  firstly, by an infinity algebraic structure   \cite{LRS} but then the problem arose that such PDEs as Euler or Navier-Stokes were not  obvious to write  in that enlarged context. A second caveat is that the homology of this complex depends on the parity of the period in each direction.
 
 \bigskip 
 The point of this paper is to  further develop this combinatorial intersection  product  to  restore the product rule for the boundary operator with  a  better approximation,  truer to  the geometric  intersection product.   The product will respect  the  intersections that are geometrically transverse and in  general position.

A parallelopiped (of any dimension) whose edges are parallel to the coordinate axes and whose vertices lie in the lattice will be called a {\it cuboid}. Such a cuboid is geometrically a Cartesian product of singletons and intervals (all of whose delimiters lie in $h\Z$). Replacing one or more of the intervals defining a cuboid by singletons at one of the endpoints of the corresponding intervals will lead to geometric objects which are cuboids of lower dimension, namely faces, edges and vertices of the original cuboid; we will call them {\it generalised faces} of the cuboid, see Figure 3. 
\[\hskip4em
\includegraphics[width=.3\textwidth]{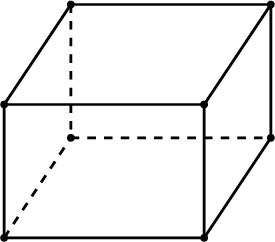}\]
\[\hbox{\small\sl Figure 3: Generalised faces of a cuboid}\]

\medskip Now consider a pair of cuboidal cells. The geometric intersection of two cuboidal cells is considered to be {\it transverse} if the set-wise intersection of the closed cells is non-empty while their tangent spaces generate the entire (three-dimensional) tangent space; for example, two intersecting lines are not transverse in three-dimensions. We say that two cuboidal cells are in {\it general position} if they intersect transversely, and in addition whenever we replace either or both cuboid by one of its generalised faces, all such pairs of cuboidal cells are either disjoint or have transverse intersection. The possible configurations of pairs of  cuboid $2h$-cells in three dimensions in general position are shown in Figure 4. 
\[\includegraphics[width=.46\textwidth]{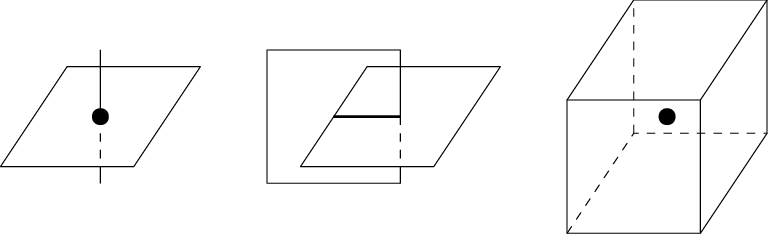}\hskip2em
\includegraphics[width=.46\textwidth]{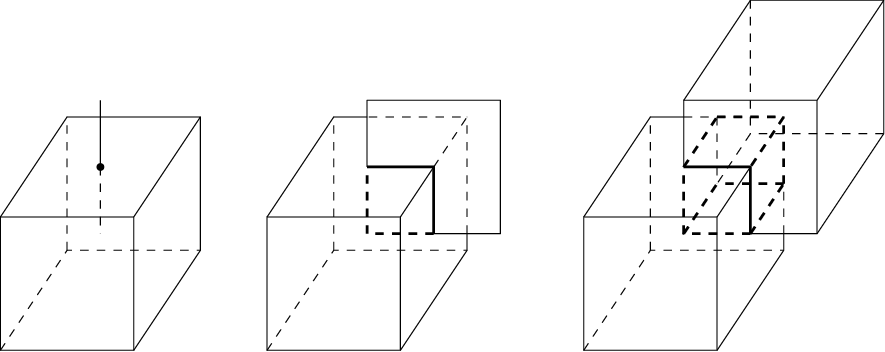}
\]
\[\hbox{\small\sl Figure 4: Intersections in the $2h$ complex in general position}\]

\medskip\noindent Observe that the geometric intersection of cuboidal cells in general position is also a cuboidal cell whose dimension is the sum of the dimensions of the initial cells minus three; equivalently the codimension of the intersection is the sum.

The possible types of configurations of $h$-cells in three dimensions which intersect transversally but are {\sl not} in general position are shown in Figure 5.
\[\includegraphics[width=.19\textwidth]{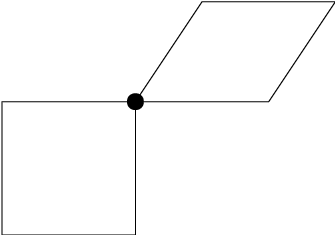}\hskip3ex
\includegraphics[width=.11\textwidth]{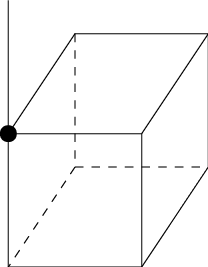}\hskip4ex
\includegraphics[width=.4\textwidth]{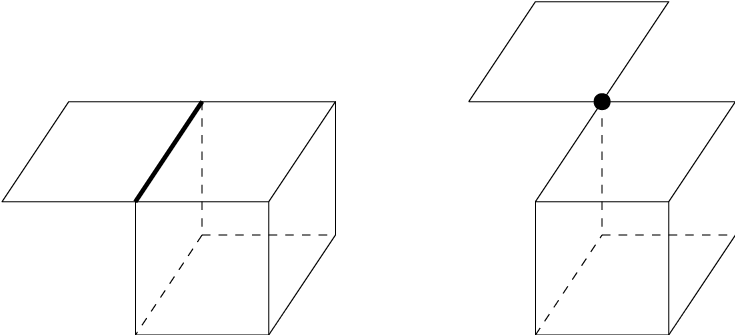}\]
\[\includegraphics[width=.78\textwidth]{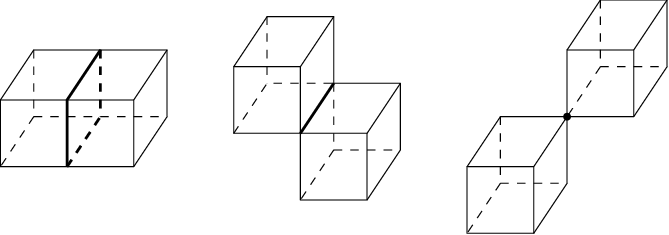}\]
\[\hbox{\small\sl Figure 5: Some cuboid intersections not in general position}\]

\medskip\noindent In the first example, the intersection `should' be of dimension $2+2-3=1$ but is geometrically a point. Intersecting cuboids of different edge lengths leads to further cases of transverse intersections not in general position. Note that any transverse intersection which is not in general position can be changed to general position by small relative translation of the cuboids. Hence the idea of this paper to wiggle the bounding faces of the cuboids by some small amount and then perform the intersections, while keeping track of the resulting probability distributions. We find that even if we start with uniform distributions for wiggling, after multiple intersections these distributions become far from uniform and indeed build an infinite hierarchy of elements.

In order to deal with this, we will work in one dimension and the resulting model can be tensored up to the needed dimension. 
The result will be an infinite algebraic structure extending the usual chain complex of a cubical lattice whose multiplication is described by geometric intersection in the case of intersecting cuboids in general position, and which is commutative, associative and satisfies  the product rule with respect to the boundary operator.

Note that the model obtained is distinct from that in \cite{ALS} where we produced a {\it finite-dimensional} transverse intersection algebra which was commutative and associative but satisfied Leibniz only on the original complex.

\begin{theorem}
The constructed transverse intersection algebra TIA is a dga over the rational numbers (that is, a graded commutative, associative algebra over the rationals satisfying Leibniz for $\d$).  TIA is generated as a linear space over the rationals by the geometric convex pieces of the original decomposition plus ideal elements (both finite in number) decorated by $2d$-tuples of non-negative integers, making it infinite dimensional.  TIA is finitely generated as an algebra over the rationals by the geometric pieces decorated by $2d$-tuples which are all zero.
\end{theorem}

{\bf Remark:A} There is a  linear chain mapping forgetting the decoration from TIA to a previous transverse intersection algebra EC (constructed in \cite{ALS}) which EC is a finite-dimensional commutative and associative algebra satisfying the product rule on all pairs of elements from the original complex, being TIA minus the decoration and the ideal elements.The EC algebra structure is the precursor of the algebra structure on TIA . This evolution was needed to improve the product rule and to  try to enable more stable fluid algebra computations. 

\

{ \bf Remark:B } Those computations based on EC showed an instability in energy even though the system was mathematically conservative. There were two likely suspects for this instability in those computations: the odd subdivision  (introduced to make the inner product  of the fluid algebra (see \cite{S10}) nondegenerate and a dangerous structure constant in the  EC algebra venturing near a pole.
 The first can be eliminated by doing even subdivisions because in TIA the inner product is essentially non degenerate for even subdivisions.  The second suspect is  buffered away from the pole in TIA.  All of this in even period decompositions for   fluid algebra computations with the TIA  discretization; and these  will be made  when  the coding of TIA is completed.

\section{The wiggling model}
Let $\Lambda$ be a $d$-dimensional lattice (periodic or infinite). Let $\epsilon>0$ be sufficiently small that $2\epsilon$ is less than the distance between any two points in $\Lambda$. 

Let $X$ be a set of (convex) polyhedral cells whose vertices lie in $\Lambda$ and which is closed under boundary and intersection in the extended sense. That is, for any $A\in{}X$, its geometric boundary $\d{}A$ is a union of non-overlapping cells in $A$ while for for any $A,B\in{}X$, either $A$ and $B$ are disjoint or their geometric intersection can be expressed as a union of non-overlapping elements of $X$.

\[\includegraphics[width=.1\textwidth]{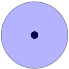}\hskip4ex
\includegraphics[width=.3\textwidth]{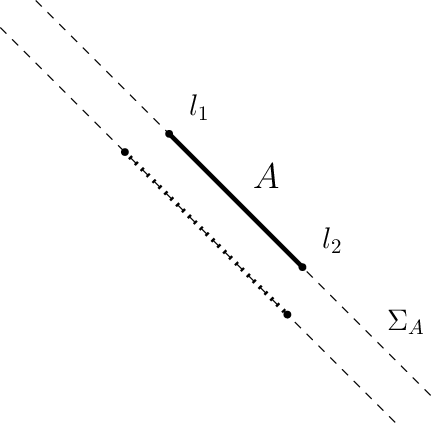}\hskip7ex
\includegraphics[width=.3\textwidth]{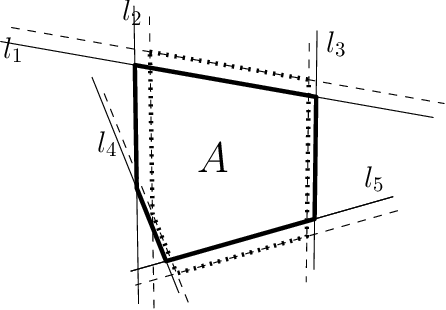}
\]
\[\hbox{\small\sl Figure 6: Wiggling a point, stick or polygon}\]

We say that $A,B\in{}X$ intersect transversally if their closures have a non-zero intersection while their normal spaces $NA,NB\subseteq\R^d$ have zero intersection.

Pick $A\in{}X$ and let $k\leq{}d$ be its dimension. Such a convex polytope $A$ can be specified by the ($k$-dimensional) affine space it generates $\Sigma_A\subseteq\R^d$ and bounding codimension-one affine subspaces $l_1,\ldots,l_N\subset\Sigma_A$.  A {\sl wiggled} version of $A$ will be the convex polytope in the affine space obtained by translating $\Sigma_A$ orthogonal to itself by at most $\epsilon$ and each of the codimension-one faces of $A\subset\Sigma_A$ similarly orthogonal to themselves (inside the new affine space) by at most $\epsilon$.  By a {\sl wiggling} of $A$ will be meant a choice of probability distribution on an $\epsilon$-ball in the normal space $NA$ plus a choice of probability distribution on the $[-\epsilon,\epsilon]$ for each (codimension one) face of $A$. The set of wiggled versions $W(A)$ of $A$ inherits a probability distribution $\mu_w$ from a wiggling $w$ of $A$ as the product distribution. By abuse of notation we will denote a wiggling of $A$ as the associated linear combination (smearing) of wiggled versions of $A$,
\[w=\int_{W(A)}B\>d\mu_w(B)\eqno{(1)}\]

\noindent The boundary of a wiggling of $A$ is defined as the linear combination of the geometric boundary of a wigged version $B$ of $A$ weighted by the probability distribution,
\[
\d(w)=\int_{W(A)}(\d{B})\>d\mu_w(B)\eqno{(2)}
\]

\noindent Given $A,A'\in{}X$ intersecting transversally and wigglings $w$ on $A$ and $w'$ on $A'$, define their transverse intersection by
\[ w\tm{}w'=\int_{W(A)\times{}W(A')} (B\cap{}B')d\mu_w(B)d\mu_{w'}(B')\eqno{(3)}\]

\noindent{\bf Remark:} Note that the result of the transverse intersection $w\tm{}w'$ may not be a (linear combination) of wigglings of elements of $X$ because the resulting probability distribution may not be a product. This typically occurs when the codimension of $A\cap{}A'$ is higher than the sum of the codimensions of $A$ and $A'$, so that the resulting object will be an `infinitesimal' object whose dimension is higher than its geometric dimension. The result of transverse intersection can however always be considered as a probability distribution on a finite-dimensional set  of shapes as in (1) and the boundaries and transverse intersections of such generalised wigglings may still be defined by (2) and (3).

\bigskip\noindent{\bf Proof of main properties}

\noindent 
Since the geometric intersection is graded commutative, and associative it follows from the definition (3) that the same is true of the transverse intersection.  Indeed a higher order transverse intersection of any number of elements $A_1,\ldots,A_m\in{}X$ of $X$ can be defined by
\[
w_1\tm\cdots\tm{}w_m=\int_{W(A_1)\times\cdots\times{}W(A_m)}(B_1\cap\cdots\cap{}B_m)d\mu_{w_1}(B_1)\cdots{}d\mu_{w_m}(B_m)
\]
Since the geometric intersection and geometric boundary satisfy Leibniz, so do $(\tm,\d)$ as defined in (2),(3),
\begin{align*}
    \d(w\tm{}w')&=\int_{W(A)\times{}W(A')} \d(B\cap{}B')d\mu_w(B)d\mu_{w'}(B')\\
&=\int_{W(A)\times{}W(A')} (\d{B}\cap{}B'+(-1)^{c_B}B\cap\d{B'})d\mu_w(B)d\mu_{w'}(B')\\
&=(\d{w})\tm{}w'+(-1)^{c_w}w\tm\d{w'}
\end{align*}

\medskip\noindent{\bf Wiggling cuboids}

\noindent A cuboid $A$ (of arbitrary dimension) as defined in the previous section, will generate an affine space$\Sigma_A$  which is a translation of one of the coordinate planes (in the generalised sense) and is delineated by pairs of codimension-one affine subspaces of $\Sigma_A$ (again parallel to coordinate axes/planes). Since all metrics are equivalent in finite dimensions, for convenience we will use the $l_\infty$ metric so that the $\epsilon$-ball becomes an $\epsilon$-cube. A wiggled version of a cuboid is another cuboid, in a parallel plane and delineated by parallel plane boundaries. A wiggling of a cuboid is determined by wigglings of all its one-dimensional projections on axes. Indeed a cuboid is the Cartesian product of singletons and intervals and a wiggled version of a cuboid is a Cartesian product of wiggled versions of these components. In particular, in one-dimension, a wiggled version of a point is another point within $\epsilon$ of the first while a wiggled version of an interval is another interval whose endpoints are within $\epsilon$ of those of the initial interval. Thus a wiggling of a point is a distribution on $[-\epsilon,\epsilon]$ while a wiggling of an interval is distribution on $[-\epsilon,\epsilon]\times[-\epsilon,\epsilon]$.

\section{Resulting one-dimensional wiggling model}
 In this section we give the result generated by repeated application of applying the procedure of the last section in one-dimension to points and intervals wiggled according to a uniform distribution. The proof is in section 4.  It is a commutative, associative, Leibniz model of the one-dimensional lattice (infinite or periodic) with parameter $h$. Denote the lattice by $\Lambda$. Choose a parameter $\epsilon>0$ so that $\epsilon<\frac{h}{2}$.

We construct a graded algebra $A$ (graded by codimension) as follows. Dimension zero objects in $A$ will be linear combinations of `extended points'; namely a basis will be given by $\{\emp_a^{m,n}|a\in\Lambda,\ m,n\in\Z^{\geq0}\}$.  The object $\emp_a^{m,n}$ should be considered as an object localized at $a\in\Lambda$ which is jiggled according to a probability distribution on $[-\epsilon,\epsilon]$ described by the parameters $m,n\in\Z^{\geq0}$,
\[f_{m,n}(z)=\frac{(m+n+1)!}{m!n!(2\epsilon)^{m+n+1}}(z+\epsilon)^m(\epsilon-z)^ne \eqno{(1)}\]
The normalization has been chosen so that $\int\limits_{-\epsilon}^{\epsilon}f_{m,n}(z)dz=1$. That is, $\emp_a^{m,n}$ can be viewed as located at the point $z$ with probability distribution $f_{m,n}(z-a)$.  Note that $f_{0,0}(z)=\frac{1}{2\epsilon}$ so that $\emp_a^{0,0}$ is a point with a uniform distribution on $[a-\epsilon,a+\epsilon]$.

\[\includegraphics[width=.8\textwidth]{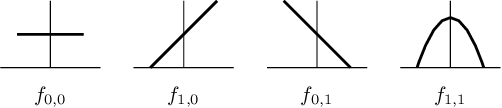}
\]
\[\hbox{\small\sl Figure 7: Point wiggling distributions}\]

\medskip
There will be two types of one-dimensional cell in $A$, namely {\sl regular jiggled intervals} and {\sl infinitesimal jiggled intervals}.  A regular jiggled interval is denoted $x_{a,b}^{m,n}$ and should be visualized as geometrically an interval $[a,b]$ on the line (with $a,b\in\Lambda$, $a<b$) in which the two end-points $a$ and $b$ are jiggled according to independent probability distributions, $f_{m,0}$ and $f_{0,n}$ respectively. That is, it is represented by the interval $[z_1,z_2]$ where the joint probability distribution of $(z_1,z_2)$ is $f_{m,0}(z_1-a)f_{0,n}(z_2-b)$.  

An infinitesimal jiggled interval is denoted $x_{a,a}^{m,n}$ and should be visualized as an infinitesimal interval around the point $a\in\Lambda$. It represents an interval $[z_1,z_2]$ with joint probability distribution $g_{m,n}(z_1-a,z_2-a)$ where
\[g_{m,n}(z_1,z_2)=\left\{\begin{array}{cl}
    \frac{(m+n+2)!}{m!n!(2\epsilon)^{m+n+2}}(z_1+\epsilon)^m(\epsilon-z_2)^n & \hbox{for $-\epsilon\leq{}z_1<z_2\leq\epsilon$} \\
     0& \hbox{otherwise}
\end{array} 
\right.\]
Again the normalization has been chosen so that $\iint\limits_{[-\epsilon,\epsilon]\times[-\epsilon,\epsilon]}g(z_1,z_2)dz_1dz_2=1$. 
 Note that in this case the distributions of the endpoints are not independent (because of the condition $z_1<z_2$). 
 \[
\includegraphics[width=.8\textwidth]{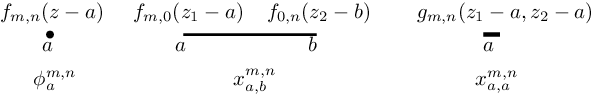}
\]
\[\hbox{\small\sl Figure 8: Generating cells}\]

\bigskip
The geometric boundary of an interval $[z_1,z_2]$ is the difference of the endpoints. Viewing the wiggled intervals as continuous linear combinations of ordinary intervals weighted by the probability distributions given, we find that the boundary map is given by
\begin{align*}
    \d(x_{ab}^{m,n})&=\mt_b^{0,n}-\mt_a^{m,0}\>,\qquad{}a<b\\
\d(x_{aa}^{m,n})&=\mt_a^{m+1,n}-\mt_a^{m,n+1}
\end{align*}
Now for intersections. Wiggled cells are almost always transverse and so we obtain a transverse intersection multiplication. The possible non-trivial intersections of zero and one-dimensional objects come in three types.
\[\includegraphics[width=.7\textwidth]{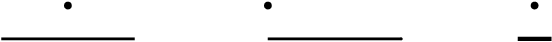}\]
\[\hbox{\small\sl Figure 9: Intersections of zero and one-dimensional objects}\]
\begin{align*}
\mt_c^{m,n}\tm{}x_{a,b}^{m',n'}&=\mt_c^{m,n}\>,\quad\hbox{if $a<c<b$}\\
\mt_a^{m,n}\tm{}x_{a,b}^{m',n'}&=\frac{(m+m'+1)!(m+n+1)!}{m!(m+n+m'+2)!}\mt_a^{m+m'+1,n}\\
\mt_b^{m,n}\tm{}x_{a,b}^{m',n'}&=\frac{(n+n'+1)!(m+n+1)!}{n!(m+n+n'+2)!}\mt_b^{m,n+n'+1}\\
\mt_a^{m,n}\tm{}x_{a,a}^{m',n'}&=\binom{m+m'+1}{m}\binom{n+n'+1}{n}\binom{m+n+m'+n'+3}{m+n+1}^{-1}\mt_a^{m+m'+1,n+n'+1}
\end{align*}

\break
For intersections of regular wiggled intervals, the possible configurations are as follows.
\[\includegraphics[width=.9\textwidth]{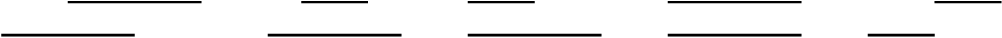}\]
\[\hbox{\small\sl Figure 11: Types of intersections of pairs of regular wiggled intervals}\]
Here is the list of intersections
\begin{align*}
x_{a,c}^{m,n}\tm{}x_{b,d}^{m',n'}&=x_{b,c}^{m',n}\>,\quad\hbox{for $a<b<c<d$}\\
x_{a,d}^{m,n}\tm{}x_{b,c}^{m',n'}&=x_{b,c}^{m',n'}\>,\quad\hbox{for $a<b<c<d$}\\
x_{a,b}^{m,n}\tm{}x_{a,c}^{m',n'}&=x_{a,b}^{m+m'+1,n}\>,\quad\hbox{for $a<b<c$}\\
x_{a,c}^{m,n}\tm{}x_{b,c}^{m',n'}&=x_{b,c}^{m',n+n'+1}\>,\quad\hbox{for $a<b<c$}\\
x_{a,b}^{m,n}\tm{}x_{a,b}^{m',n'}&=x_{a,b}^{m+m'+1,n+n'+1}\\
x_{a,b}^{m,n}\tm{}x_{b,c}^{m',n'}&=\frac{(m'+1)!(n+1)!}{(m'+n+2)!}x_{b,b}^{m',n}\>,\quad\hbox{for $a<b<c$}
\end{align*}

Finally we have intersections of one-dimensional objects which involve infinitesimals
\[\includegraphics[width=.6\textwidth]{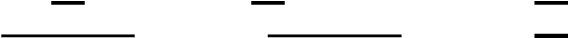}\]
\[\hbox{\small\sl Figure 12: Types of intersections of involving infinitesimals}\]
The values of these intersections are
\begin{align*}
x_{b,b}^{m,n}\tm{}x_{a,c}^{m',n'}&=x_{b,b}^{m,n}\>,\quad\hbox{for $a<b<c$}\\
x_{a,a}^{m,n}\tm{}x_{a,b}^{k,l}&=\frac{(m+n+2)!(m+k+2)!}{(m+1)!(m+n+k+3)!}x_{a,a}^{m+k+1,n}\\
x_{b,b}^{m,n}\tm{}x_{a,b}^{k,l}&=\frac{(m+n+2)!(n+l+2)!}{(n+1)!(m+n+l+3)!}x_{b,b}^{m,n+l+1}\\
x_{a,a}^{m,n}\tm{}x_{a,a}^{m',n'}&=\binom{m+n+2}{m+1}\binom{m'+n'+2}{m'+1}\binom{m+n+m'+n'+4}{m+m'+2}^{-1}x_{a,a}^{m+m'+1,n+n'+1}
\end{align*}

\break\section{Derivation of the one-dimensional model}
In this section we derive the one-dimensional given explicitly in the last section from the general strategy of probabilistic wiggling described in \S2 applied to a one-dimensional lattice.

In a one-dimensional lattice $\Lambda$ we start with basic objects which are points $a$ ($a\in\Lambda$) and intervals $[a,b]$ with $a,b\in\Lambda$. A wiggled version of the point $a$ is a (random) point in the interval $[a-\epsilon,a+\epsilon]$ and a wiggled version of the interval $[a,b]$ is an interval $[z_1,z_2]$ where $z_1\in[a-\epsilon,a+\epsilon]$ and $z_2\in[b-\epsilon,b+\epsilon]$.

Any probability distribution in $[a-\epsilon,a+\epsilon]$ will define a wiggling of $a$, that is, a continuous linear combination of points near $a$. 

For $m,n\in\Z^{\geq0}$, consider the wiggling of $a$ defined by a probability distribution on $[a-\epsilon,a+\epsilon]$ which is given by a polynomial in $z$ vanishing to order $m$ at $z=a-\epsilon$ and to order $n$ at $z=a+\epsilon$, that is proportional to $(z-a+\epsilon)^m(a+\epsilon-z)^n$. Since 
\[\int_{-\epsilon}^\epsilon(z-\epsilon)^m(\epsilon-z)^ndz=\epsilon^{m+n+1}\int_{-1}^1(1+z)^m(1-z)^ndz=\frac{m!n!(2\epsilon)^{m+n+1}}{(m+n+1)!}\]
thus the appropriate probability distribution is
\[f_{m,n}(z-a)=\frac{(m+n+1)!}{m!n!(2\epsilon)^{m+n+1}}(z-a+\epsilon)^m(\epsilon+a-z)^n\]

\noindent{\bf Definition:} Denote by $\mt_a^{m,n}$ the wiggling of the point $a$ which is specified as the linear combination of points $z$ with distribution $f_{m,n}(z-a)$.

Define a wiggling of the interval $[a,b]$, in which the joint probability distribution of the endpoints $z_1,z_2$ of the interval $[z_1,z_2]$ is given by 
\[
f_{m,0}(z_1-a)f_{0,n}(z_2-b)=\frac{(m+1)(n+1)}{(2\epsilon)^{m+n+2}}(z_1-a+\epsilon)^m(\epsilon+b-z_2)^n
\]
\noindent{\bf Definition:} Denote by $x_{a,b}^{m,n}$ the wiggling of the interval $[a,b]$ which is specified as the linear combination of intervals $[z_1,z_2]$ with joint distribution $f_{m,0}(z_1-a)f_{0,n}(z_2-b)$.

\medskip
\noindent{\bf Intersection of points and intervals $\mt_c^{m,n}\tm{}x_{a,b}^{m',n'}$}

\noindent Consider the transverse intersection $\mt_c^{m,n}\tm{}x_{a,b}^{m',n'}$. We take the intersection of wiggled versions 
\[
\{z\}\cap[z_1,z_2]=\left\{\begin{array}{cc}
     \{z\}&\hbox{if $z_1\leq{}z\leq{}z_2$}  \\
     \mt& \hbox{otherwise}
\end{array}
\right. 
\]
in which the random variables $z,z_1,z_2$ have joint probability distribution $f_{m,n}(z-c)f_{m',0}(z_1-a)f_{0,n'}(z_2-b)$. Thus we obtain as intersection the point $z$ with probability distribution
\begin{align*}
   & \int_{-\infty}^z\int_z^\infty{}f_{m,n}(z-c)f_{m',0}(z_1-a)f_{0,n'}(z_2-b)dz_2dz_1\\
    &=f_{m,n}(z-c)\left(\int_{-\infty}^zf_{m',0}(z_1-a)dz_1\right)\left(\int_z^\infty{}f_{0,n'}(z_2-b)dz_2\right)
\end{align*}
Recall that $a,b,c\in\Lambda$ and the lattice spacing is such that $2\epsilon<h$ so that the only configurations are those of general position $c<a$, $a<c<b$, $c>b$ and the cases of coincident points $c=a$  or $c=b$. The evaluations here are
\begin{align*}
    c<a\hbox{ or }c>b&:\quad0\\
    a<c<b&:\quad{}f_{m,n}(z-c)\\
    c=a&:\quad{}f_{m,n}(z-a)\left(\int_{a-\epsilon}^zf_{m',0}(z_1-a)dz_1\right)\\
    c=b&:\quad{}f_{m,n}(z-a)\left(\int_z^{a+\epsilon}{}f_{0,n'}(z_2-b)dz_2\right)
\end{align*}
In the case $c=a$, We find that
\[
\int_{a-\epsilon}^zf_{m',0}(z_1-a)dz_1
=(m'+1)(2\epsilon)^{-m'-1}\int_{-\epsilon}^{z-a}(u+\epsilon)^{m'}du
=(2\epsilon)^{-m'-1}(z-a+\epsilon)^{m'+1}
\]
so that 
\begin{align*}
  f_{m,n}(z-a)\left(\int_{a-\epsilon}^zf_{m',0}(z_1-a)dz_1\right)
&=\frac{(m+n+1)!}{m!n!(2\epsilon)^{m+m'+n+2}}(z-a+\epsilon)^{m+m'+1}(a+\epsilon-z)^n \\
&=\frac{(m+m'+1)!(m+n+1)!}{m!(m+m'+n+2)!}f_{m+m'+1}{n}(z-a)
\end{align*}
Similarly for the case $c=b$ so that we finally obtain
\[
\mt_c^{m,n}\tm{}x_{a,b}^{m',n'}=\left\{\begin{array}{ll}
     0&\hbox{if $c<a$ or $c>b$}  \\
     \mt_c^{m,n}& \hbox{if $a<c<b$}\\
     \frac{(m+m'+1)!(m+n+1)!}{m!(m+m'+n+2)!}\mt_a^{m+m'+1,n}&\hbox{if $c=a$}\\
     \frac{(n+n'+1)!(m+n+1)!}{n!(m+n+n'+2)!}\mt_a^{m,n+n'+1}&\hbox{if $c=b$}
\end{array}
\right. 
\]

\bigskip
\noindent{\bf Intersection of pairs of intervals}

\noindent Consider an intersection of a pair of intervals $x_{a,b}^{m,n}\tm{}x_{a',b'}^{m',n'}$.  This is defined by taking wiggled versions $[z_1,z_2]$ and $[z'_1,z'_2]$ of the two intervals and intersecting them, then taking their linear combination according to the joint probability distribution of $z_1,z_2,z'_1,z'_2$ (all independent). As in the previous calculation, when the initial (unwiggled) intersection is either empty or in general position, the computation is immediate. Intersection is commutative leaving three cases of this sort
\[\includegraphics[width=.7\textwidth]{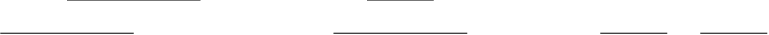}\]
giving intersections
\[
x_{a,b}^{m,n}\tm{}x_{a',b'}^{m',n'}=\left\{\begin{array}{ll}
     x_{a',b}^{m',n}&\hbox{if $a<a'<b<b'$}  \\
     x_{a',b'}^{m',n'}& \hbox{if $a<a'<b'<b$}\\
     0&\hbox{if $[a,b]\cap[a',b']=\mt$}
\end{array}
\right. 
\]
The cases of intersections which are not in general position are those in which the intervals share one or both endpoints,
\[\includegraphics[width=.9\textwidth]{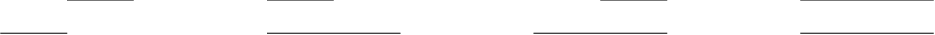}\]
(i) For the first case, $x_{a,b}^{m,n}\tm{}x_{b,c}^{m',n'}$, for $a<b<c$,  we see that the wiggled intervals $[z_1,z_2]$ of $[a,b]$ and $[z'_1,z'_2]$ of $[b,c]$ intersect only when $z'_1<z_2$ in which case their intersection is $[z'_1,z_2]$. The joint probability distribution describing the frequency with which this particular interval occurs is now given by integrating over the other two variables $z_1$, $z'_2$
\[
\iint{}f_{m,0}(z_1-a)f_{0,n}(z_2-b)f_{m',0}(z'_1-b)f_{0,n'}(z'_2-c)dz_1dz'_2=f_{m',0}(z'_1-b)f_{0,n}(z_2-b)
\]
so that we obtain an interval around $b$ with weighting as given. Changing the names of the variables, the result of the intersection is the interval $[z_1,z_2]$ with distribution
\[
\left\{
\begin{array}{cc}
(m'+1)(n+1)(2\epsilon)^{-m'-n-2}(z_1-b+\epsilon)^{m'}(b+\epsilon-z_2)^n&\hbox{if $z_1<z_2$}\\
0&\hbox{if $z_1>z_2$}
\end{array}\right.
\]
The total probability of a non-empty result is less than one. We define an infinitesimal wiggled interval around a point $a$ by normalising such a distribution. Calculating
\[\int_{-\epsilon}^{\epsilon}\int_{z_1}^\epsilon
(z_1+\epsilon)^m(\epsilon-z_2)^ndz_2dz_1=\frac{m!n!}{(m+n+2)!}(2\epsilon)^{m+n+2}\]
we deduce the correct normalisation.

\noindent{\bf Definition:} Denote by $x_{a,a}^{m,n}$ the wiggling of the infinitesimal interval around $a$ which is specified as the linear combination of intervals $[z_1,z_2]$ with joint distribution $g_{m,n}(z_1-a,z_2-a)$ where 
\[
g_{m,n}(z,w)=\left\{\begin{array}{cl}
    \frac{(m+n+2)!}{m!n!(2\epsilon)^{m+n+2}}(z+\epsilon)^m(\epsilon-w)^n & \hbox{for $-\epsilon\leq{}z<w\leq\epsilon$} \\
     0& \hbox{otherwise}
\end{array} 
\right.
\]
We conclude that 
\[
x_{a,b}^{m,n}\tm{}x_{b,c}^{m',n'}=\frac{(m'+1)!(n+1)!}{(m'+n+2)!}x^{m',n}_{b,b}=\binom{m'+n+2}{n+1}^{-1}x^{m',n}_{b,b}
\]

\medskip\noindent (ii) For the second case, we have an intersection $x_{a,b}^{m,n}\tm{}x_{a,c}^{m',n'}$ with $a<b<c$. Wiggled versions of the intervals $[a,b]$ and $[a,c]$ will be $[z_1,z_2]$ and $z'_1,z'_2]$ with $z_1,z'_1<z_2<z'_2$ and therefore their intersection will be $[\max(z_1,z'_1),z_2]$. This will be $[z,w]$ in either of the two cases $z_1<z'_1=z<z_2=w<z'_2$ or $z'_1<z_1=z<z_2=w<z'_2$. That is the result is $[z,w]$ with joint probability distribution non-zero for $|z-a|\leq\epsilon$, $|w-b|\leq\epsilon$,
\begin{align*}
 &\int_{a-\epsilon}^z\int_{c-\epsilon}^{c+\epsilon}f_{m,0}(z_1-a)f_{0,n}(w-b)f_{m',0}(z-a)f_{0,n'}(z'_2-c)dz'_2dz_1  \\
 &+\int_{a-\epsilon}^z\int_{c-\epsilon}^{c+\epsilon}f_{m,0}(z-a)f_{0,n}(w-b)f_{m',0}(z'_1-a)f_{0,n'}(z'_2-c)dz'_2dz'_1  \\
 &=f_{0,n}(w-b)f_{m',0}(z-a)\left(\int_{a-\epsilon}^zf_{m,0}(z_1-a)dz_1\right)\\
 &\qquad+f_{m,0}(z-a)f_{0,n}(w-b)\left(\int_{a-\epsilon}^zf_{m',0}(z'_1-a)dz'_1\right)\\
 &=f_{m+m'+1,0}(z-a)f_{0,n}(w-b)
\end{align*}
where in the last step we use that $\int_{a-\epsilon}^zf_{m,0}(z_1-a)dz_1=(2\epsilon)^{-m-1}(z-a)^{m+1}$. The conclusion is that
\[
x_{a,b}^{m,n}\tm{}x_{a,c}^{m',n'}=x_{a,c}^{m+m'+1,n}\hbox{ for }a<b<c
\]
(iii) Similarly we have the third case
\[
x_{a,c}^{m,n}\tm{}x_{b,c}^{m',n'}=x_{b,c}^{m',n+n'+1}\hbox{ for }a<b<c
\]

\noindent(iv) In the case of an interval intersected with itself $x_{a,b}^{m,n}\tm{}x_{a,b}^{m',n'}$, a pair of wiggled versions $[z_1,z_2]$ and $[z'_1,z'_2]$ of $[a,b]$ will intersect in $[\max(z_1,z'_1),\min(z_2,z'_2)]$. This will be the interval $[z,w]$ in one of four cases with either $z_1<z=z'_1$ or $z'_1<z=z_1$ and either $w=z_2<z'_2$ or $w=z'_2<z_2$. The result of the intersection thus will be $[z,w]$ with joint probability distribution
\begin{align*}
   & \left[f_{m',0}(z-a)\left(\int_{a-\epsilon}^zf_{m,0}(z_1-a)dz_1\right)+f_{m,0}(z-a)\left(\int_{a-\epsilon}^zf_{m',0}(z'_1-a)dz'_1\right)\right]\\
  &\cdot\left[f_{0,n}(w-b)\left(\int_w^{b+\epsilon}f_{0,n'}(z'_2-b)dz'_2\right)+f_{0,n'}(w-b)\left(\int_w^{b+\epsilon}f_{0,n}(z_2-b)dz_2\right)\right]\\
  &=f_{m+m'+1,0}(z-a)f_{0,n+n'+1}(w-b)
\end{align*}
from which we derive the result
\[
x_{a,b}^{m,n}\tm{}x_{a,b}^{m',n'}=x_{a,b}^{m+m'+1,n+n'+1}
\]

\bigskip
\noindent{\bf Intersections involving infinitesimal intervals}

\noindent Now that we introduced a new object, an infinitesimal interval $x_{a,a}^{m,n}$ we need to discuss transverse intersections between it and other objects, points, intervals and other infinitesimal intervals.
\[\includegraphics[width=.9\textwidth]{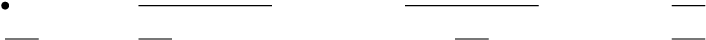}\]
(i) To intersect a point and get a non-trivial result, the infinitesimal interval must be located at the same point $\mt_a^{m,n}\tm{}x_{a,a}^{m',n'}$. A wiggled version of the point $z$ and of the infinitesimal interval $[z_1,z_2]$ will intersect precisely when $z_1\leq{}z\leq{}z_2$ and then the intersection will be the point $z$. So the probability associated with $z$ is
\begin{align*}
    &\int_{a-\epsilon}^z\int_z^{a+\epsilon}f_{m,n}(z-a)g_{m',n'}(z_1-a,z_2-a)dz_2dz_1\\
&=\frac{(m+n+1)!(m'+n'+2)!}{m!n!m'!n'!(2\epsilon)^{m+n+m'+n'+3}}(z-a-\epsilon)^m(a+\epsilon-z)^n\\
&\qquad\int_{-\epsilon}^{z-a}\int_{z-a}^\epsilon (z_1-a+\epsilon)^{m'}(\epsilon+a-z_2)^{n'}dz_2dz_1\\
&=\frac{(m+n+1)!(m'+n'+2)!}{m!n!(m'+1)!(n'+1)!(2\epsilon)^{m+n+m'+n'+3}}(z-a-\epsilon)^{m+m'+1}(a+\epsilon-z)^{n+n'+1}\\
&=\binom{m+m'+1}{m'}\binom{n+n'+1}{n'}\binom{m+n+m'+n'+3}{m'+n'+1}^{-1}f_{m+m'+1,n+n'+1}(z)
\end{align*}
from which we conclude that
\[
\mt_a^{m,n}\tm{}x_{a,a}^{m',n'}=\binom{m+m'+1}{m'}\binom{n+n'+1}{n'}\binom{m+n+m'+n'+3}{m'+n'+1}^{-1}\mt_a^{m+m'+1,n+n'+1}
\]


\medskip\noindent (iii) An intersection of form $x_{b,b}^{m,n}\tm{}x_{a,c}^{m',n'}$ for $a<b<c$ is in general position and therefore immediately reduces to $x_{b,b}^{m,n}$.

\noindent The derivation of the formulae for (ii) and (iv) are similar.

\bigskip\noindent{\bf Boundaries}

\noindent For the last part of the data, we compute the boundaries of wiggled intervals and of infinitesimal wiggled intervals. To compute $\d(x_{a,b}^{m,n})$, we take a wiggled version $[z_1,z_2]$ of the interval, whose boundary is the difference of points $\mt_{z_2}-\mt_{z_1}$. This is to be weighted by the probability distribution $f_{m,0}(z_1-a)f_{0,n}(z_2-b)$ and so 
\begin{align*}
    \d(x_{a,b}^{m,n})&=\int_{a-\epsilon}^{a+\epsilon}\int_{b-\epsilon}^{b+\epsilon}
f_{m,0}(z_1-a)f_{0,n}(z_2-b)(\mt_{z_2}-\mt_{z_1})dz_2dz_1\\
&=\int_{b-\epsilon}^{b+\epsilon}f_{0,n}(z_2-b)\mt_{z_2}dz_2-\int_{a-\epsilon}^{a+\epsilon}f_{m,0}(z_1-a)dz_1=\mt_b^{0,n}-\mt_a^{m,0}
\end{align*}

\noindent On the other hand, the computation of the boundary $\d(x_{a,a}^{m,n})$ of the wiggled infinitesimal interval works similarly 
\begin{align*}
    \d(x_{a,a}^{m,n})&=\int_{a-\epsilon}^{a+\epsilon}\int_{z_1}^{a+\epsilon}
g_{m,n}(z_1-a,z_2-a)(\mt_{z_2}-\mt_{z_1})dz_2dz_1\\
&=\int_{-\epsilon}^{\epsilon}\left(\int_{-\epsilon}^{z_2}g_{m,n}(z_1,z_2)dz_1\right)\mt_{a+z_2}dz_2-\int_{-\epsilon}^{\epsilon}\left(\int_{z_1}^{\epsilon}g_{m,n}(z_1,z_2)dz_2\right)\mt_{a+z_1}dz_1
\end{align*}
Observe that
\begin{align*}
    \int_{-\epsilon}^{z_2}g_{m,n}(z_1,z_2)dz_1&=\frac{(m+n+2)!}{(m+1)!n!}(2\epsilon)^{-m-n-2}(\epsilon+z_2)^{m+1}(\epsilon-z_2)^n\\
    \int_{z_1}^{\epsilon}g_{m,n}(z_1,z_2)dz_2
    &=\frac{(m+n+2)!}{m!(n+1)!}(2\epsilon)^{-m-n-2}(\epsilon+z_1)^m(\epsilon-z_1)^{n+1}
\end{align*}
It follows that $\d(x_{a,a}^{m,n})=\mt_a^{m+1,n}-\mt_a^{m,n+1}$.

\section{Connections with other work}

Here are some examples of subspaces of TIA and relations to other work.

\medskip
\noindent{\bf Example 1:} For any $K\in\N$, the subspace of TIA generated by those generators whose decorations are all at least $K$, forms an ideal with respect to transverse multiplication. It is however not closed under the boundary operator $\d$.
\medskip

\noindent{\bf Example 2:} Consider a one-dimensional periodic lattice with lattice parameter $h$. Inside the TIA defined above, we can consider the linear space $W$ spanned by decorated points and decorated intervals of length $2h$. There is an involution $*:W\to{}W$ defined by
\begin{align*}
    *(\mt_a^{m,n})&=x_{a-h,a+h}^{m,n}\>,\\
    *(x_{a-h,a+h}^{m,n})&=\mt_a^{m,n}\>.
\end{align*}
Let $W_0$ denote that part of $W$ with decorations $0,0$, that is, spanned by $\mt_a^{0,0}$ and $x_{a-h,a+h}^{0,0}$ for $a\in\Lambda$. The star operator is also an involution on $W_0$. The $k$-fold intersection of a $2h$-interval decorated by $0,0$ with itself is the same $2h$-interval decorated by $k,k$.  Hence the subalgebra $U$ of TIA generated by $W_0$ has basis $\mt_a^{m,n}$, $x_{a,a+h}^{m,n}$, $x_{a-h,a+h}^{m,m}$ for $a\in\Lambda$, $m,n\geq0$.

Taking the tensor product of three copies of $W$ yields the subspace $W^{\otimes3}$ of the three-dimensional TIA generated by $2h$-cubes (of all dimensions 0,1,2,3) on which again there is a natural involution $*$. The star operator is also an involution on $W_0^{\otimes3}$. The subalgebra of TIA generated by $W_0^{\otimes3}$ is $U^{\otimes3}$.

\medskip
\noindent{\bf Remark:} The construction of the dga TIA here can be compared and contrasted to Whitney forms \cite{W} on simplicial complexes defined by polynomial forms with $\Q$-coefficients. (This was used in \cite{S77} and earlier over the reals by Ren\'e Thom \cite{Thom2} to study Postnikov systems.) 

The dga TIA makes sense for certain cubical complexes not for arbitrary simplicial complexes where Whitney can be defined \cite{W}. Also vice versa, the idea of Whitney forms uses properties of simplices and are not easily developed for cubical complexes. Finally, Whitney forms, as in \cite{S77}, consist of  all polynomials  in several variables with $\Q$ coefficients. The dga TIA  uses particular distributions described by  specific polynomials indexed by tuples of nonnegative integers related to the wiggling. These extra parameters arise because x is almost never transverse to itself and requires  wiggling  creating these parameters.

\bigskip
\noindent{\bf Fluid algebras}

\noindent
In \cite{S10}, Sullivan reformulated Euler's fluid equation in terms of the fluid algebra of differential forms. A finite-dimensional version of the fluid equation is generated by a finite-dimensional fluid algebra and we can use a transverse intersection algebra in place of differential forms to generate such finite things.
More precisely, a {\sl fluid algebra} \cite{S10} is a vector space $V$ along with
\begin{enumerate}
    \item[1.] a positive definite inner product $(\ ,\ )$ (the {\sl metric})
    \item[2.] a symmetric non-degenerate bilinear form $\langle\ ,\ \rangle$ (the {\sl linking form})
    \item[3.] an alternating trilinear form $\{\ ,\ ,\ \}$ (the {\sl triple intersection form})
\end{enumerate}
Given a fluid algebra, the {\sl associated Euler equation} is an evolution equation for $X(t)\in{}V$ given implicitly by
\[
(\dot{X},Z)=\{X,DX,Z\}\hbox{ for all test vectors }Z\in{}V
\]
where $D:V\to{}V$ is the operator defined by $\langle{}X,Y\rangle=(DX,Y)$ for all $X,Y\in{}V$.

\medskip
\noindent{\bf Example 3:} Consider a three-dimensional periodic lattice with lattice parameter $h$ and let $V$ be the subspace of TIA consisting of coexact linear combinations of $2h$-squares in the lattice with the cells decorated by a six-tuple of zeroes. Observe that both the star and boundary operators on cells decorated by zeroes are decorated by zeroes.
Define a pre-fluid algebra on $V$ by
\begin{align*}
 (a,b)&=\#(*a\tm{}b)\>,\\
\langle{}a,b\rangle&=\#(a\tm\d{b})\>,\\
\{a,b,c\}&=\#(a\tm{}b\tm{}c)
\end{align*}
where $*$ is the star defined in Example 2 and $\#:X\to\Q$ is an augmentation, a linear map on the subspace of TIA generated by codimension three objects (that is, points). That the triple form is alternating follows from the fact that TIA is graded commutative. That the linking form is symmetric follows from the product rule in TIA so long as for all intersections $x$ of pairs of elements of $V$,
\[\#(\d(x))=0\]
This can be verified for both the standard augmentation which counts points and the modified one which takes a decorated point to a power of a parameter $\delta\in(0,1]$ given by the sum of the six decorating integers of the point.  Furthermore, the inner product $(\>,\>)$ will be symmetric. It will be positive definite in either of the two cases, odd lattice size and $\delta\leq1$ or even lattice size and $\delta<1$.

The `structure constants' entering the fluid algebra so generated come from intersections of triples of zero decorated $2h$-squares and of a pair of a $2h$ square and a $2h$-stick.  Since $2h$-squares can be considered as the Cartesian product of a point and two intervals of length $2h$, such intersections are Cartesian products of intersections of pairs or triples of elements of the one-dimensional TIA with zero decorations.  Indeed, intersections of pairs of zero decorated elements from the one-dimensional lattice $\mt_a^{0,0}$ and $x_{a,b}^{0,0}$ involve points, decorated with $0,1$ and $1,0$ 
\[\mt_a^{0,0}\tm{}x_{a,b}^{0,0}=\frac12\mt_a^{1,0}\>,\quad
    \mt_b^{0,0}\tm{}x_{a,b}^{0,0}=\frac12\mt_b^{0,1}\]
as well as additional sticks
$x_{a,b}^{1,0}, x_{a,b}^{0,1}, x_{a,b}^{1,1}$ and infinitesimal sticks $x_{a,a}^{0,0}$ from
\[ x_{a,b}^{0,0}\tm{}x_{a,c}^{0,0}=x_{a,b}^{1,0}\>,\ \
    x_{a,c}^{0,0}\tm{}x_{b,c}^{0,0}=x_{b,c}^{0,1}\>,\ \
    x_{a,b}^{0,0}\tm{}x_{a,b}^{0,0}=x_{a,b}^{1,1}\>,\ \
    x_{a,b}^{0,0}\tm{}x_{b,c}^{0,0}=\frac12x_{b,b}^{0,0}
\]
Note that even the $(0,0)$ infinitesimal stick has non-trivial boundary,
\[\d(x_{a,a}^{0,0})=\mt_a^{1,0}-\mt_a^{0,1}\]
Triple intersections of two sticks and a point from the original complex involve the additional intersections
\begin{align*}
    \mt_a^{1,0}\tm{}x_{a,b}^{0,0}&=\frac23\mt_a^{2,0},\
    \mt_a^{0,1}\tm{}x_{a,b}^{0,0}=\frac13\mt_a^{1,1},\
    \mt_b^{1,0}\tm{}x_{a,b}^{0,0}=\frac13\mt_a^{1,1},\\
    \mt_b^{0,1}\tm{}x_{a,b}^{0,0}&=\frac23\mt_a^{0,2},\
    \mt_a^{0,0}\tm{}x_{a,b}^{1,n}=\frac13\mt_a^{2,0},\
    \mt_a^{0,0}\tm{}x_{a,a}^{0,0}=\frac13\mt_a^{1,1}
\end{align*}
In this setting, the augmentation is given by $\#(\mt_a^{m,n})=\delta^{m+n}$.  These are all special cases of the formulae in the previous section.
\[\includegraphics[width=.5\textwidth]{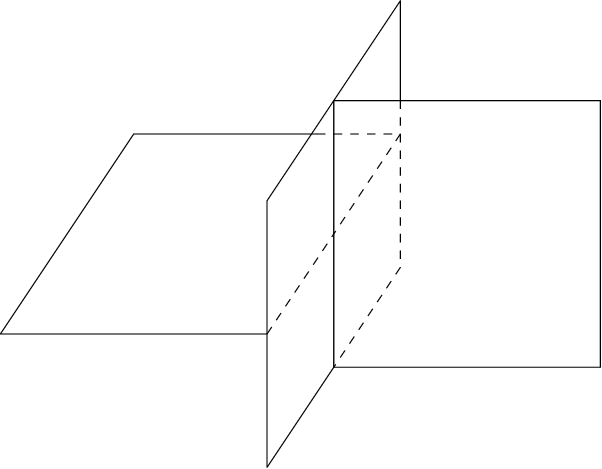}\]
As an example of a triple intersection we give
\begin{align*}
   & (\mt_0^{0,0}\otimes{}y_{-h,h}^{0,0}\otimes{}z_{-h,h}^{0,0})
\tm(x_{0,2h}^{0,0}\otimes\mt_0^{0,0}\otimes{}z_{-h,h}^{0,0})
\tm(x_{-2h,0}^{0,0}\otimes{}y_{-h,h}^{0,0}\otimes\mt_{0}^{0,0})\\
&\qquad=-\frac16\mt_0^{1,1}\otimes\mt_0^{1,1}\otimes\mt_0^{1,1}
\end{align*}
as opposed to the geometrically similar triple intersection in \cite{ALS} which had a coefficient of $-\frac14$.

\end{document}